\documentclass[11pt]{article}
\usepackage{amssymb}
\newtheorem{theorem}{Theorem}

\newtheorem{corollary}[theorem]{Corollary}

\def\qed{\hfill $\Box$\medskip}

\def\IC{{\bf C}}

\def\cH{{\cal H}}
\def\cB{{\cal B}}
\def\BH{{\cB(\cH)}}
\def\Re{{\rm Re}\,}

\def\conv{{\rm conv}\,}
\def\dim{{\rm dim}\,}

\begin{document}
\begin{center}
{\large{\bf Condition for the higher rank numerical range to be non-empty}}

\medskip
{\sc Chi-Kwong Li},\footnote{
Department of Mathematics, The College of William and Mary,
Williamsburg, VA 23185 (ckli@math.wm.edu).
Research of Li was partially supported by a USA NSF grant and
a HK RGC grant. He is an honorary professor of the University of 
Hong Kong.}
{\sc Yiu-Tung Poon}\footnote{
Department of Mathematics,
Iowa State University, 
Ames, IA 50051 (ytpoon@iastate.edu).}
{\sc and 
Nung-Sing Sze}
\footnote
{Department of Mathematics,
University of Connecticut, 
Storrs, CT 06269  (sze@math.uconn.edu).}

\date{}
\end{center}

\bigskip
\noindent
{\large{\bf Abstract}}

\sl
It is shown that 
the rank-$k$ numerical range of every $n$-by-$n$ complex matrix is
non-empty if $k < n/3+1$.  
The proof is based on a recent characterization of
the rank-$k$ numerical range by Li and Sze, the Helly's theorem
on compact convex sets, and some eigenvalue inequalities.
In particular, the result implies that $\Lambda_2(A)$ is
non-empty if $n \ge 4$. This confirms a conjecture
of Choi et al.  If $k \ge n/3+1$, 
an $n$-by-$n$ complex matrix is given for which
the rank-$k$ numerical range is empty.
Extension of the result to bounded linear operators acting on 
an infinite dimensional Hilbert space is also discussed.

\rm



\section{Introduction}

Let $M_n$ be the algebra of $n\times n$ complex matrices.
In \cite{Cet1}, the authors introduced the notion of the
{\it rank-$k$ numerical range} of
$A \in M_n$ defined and denoted by
$$\Lambda_k(A) = \{ \lambda  \in \IC: X^*AX = \lambda  I_k, \
X \hbox{ is } n\times k \hbox{ such that }  X^*X = I_k\}$$ in
connection to the study of quantum error correction; see
\cite{Cet2}. 
Evidently, $\lambda \in \Lambda_k(A)$ if and only if
there is a unitary matrix $U \in M_n$ such that $U^*AU$ has $\lambda
I_k$ as the leading principal submatrix. When $k = 1$, this concept
reduces to the classical numerical range, which is well known to be
convex by the Toeplitz-Hausdorff theorem; for example, see \cite{L}
for a simple proof. In \cite{Cet} the authors conjectured that
$\Lambda_k(A)$ is convex, and reduced the convexity problem to the
problem of showing that $0 \in \Lambda_k(A)$ for
$$A = \pmatrix{I_k & X \cr Y & -I_k \cr}$$
for arbitrary $X, Y \in M_k.$
They further reduced this problem to the existence of a Hermitian
matrix $H$ satisfying  the matrix equation
\begin{equation} \label{eq1}
I_k + MH + HM^* - HPH = H
\end{equation}
for arbitrary $M \in M_k$ and positive definite $P \in M_k$.
In \cite{W}, the author observed that equation (\ref{eq1})
can be rewritten as the continuous Riccati equation
\begin{equation} \label{eq2}
HPH-H(M^*-I_k/2)-(M-I_k/2)H - I_k = 0_k,
\end{equation}
and existing results on Riccati equation
will ensure its solvability;
for example, see \cite[Theorem 4]{LR}.
This establishes the convexity of $\Lambda_k(A)$.

For a Hermitian $X \in M_n$, let $\lambda_1(X) \ge \cdots \ge
\lambda_n(X)$ be the eigenvalues of $X$. 
In \cite{LS}, it was shown that 
\begin{equation}\label{lambda}\Lambda_k(A) = \{ \mu \in \IC: 
e^{it}\mu+ e^{-it}\overline{\mu} \le \lambda_k(e^{it}A + e^{-it}A^*)
\mbox{ for all }t\in [0,2\pi)\}\,.\end{equation} 
In particular, $\Lambda_k(A)$
is the intersection of closed half planes on $\IC$, and therefore is
always convex. Moreover, if $A \in M_n$ is normal with eigenvalues
$\lambda_1, \dots, \lambda_n$, then
$$\Lambda_k(A)
= \bigcap_{1 \le j_1 < \cdots < j_{n-k+1} \le n} \conv\{\lambda_{j_1},\dots,
\lambda_{j_{n-k+1}}\}.$$
This confirms a conjecture in \cite{Cet0}.

While many interesting results have been obtained for $\Lambda_k(A)$,
see \cite{Cet,Cet0,Cet1,Cet2} for example,
there are some basic questions whose answers are unknown.
The purpose of this paper is to answer the following.

\medskip\noindent
{\bf Problem} Determine $n$ and $k$ such that 
$\Lambda_k(A)$ is non-empty for every $A \in M_n$.

\medskip
It is well-known that the classical numerical range $\Lambda_1(A)$
is non-empty.
For $k > n/2$ , $\Lambda_k(A)$ has at most one element and one can 
easily construct $A \in M_n$ such that $\Lambda_k(A) = \emptyset$;
see Proposition 2.2 and Corollary 2.3 in \cite{Cet1}.
The situation for $\Lambda_k(A)$ with $n/2 \ge k > 1$ is not so
clear.  In \cite{Cet0}, the authors conjectured that
$\Lambda_2(A) \ne \emptyset$ for $n \ge 4$.

In the next section, we show that 
$\Lambda_k(A)$ is non-empty for 
every $A \in M_n$ if and only if $k < n/3 + 1$.
In particular, it confirms the conjecture
in \cite{Cet0} that  $\Lambda_2(A) \ne
\emptyset$ if $n \ge 4$. 
We also consider extension 
of the result to infinite dimensional bounded linear
operators.

\section{Results and proofs}

\begin{theorem}
Let $A \in M_n$, and let
$k$ be a positive integer such that $k < n/3+1$.
Then $\Lambda_k(A)$ is non-empty.
\end{theorem}

\it Proof. \rm
Evidently, $\Lambda_k(A) \subseteq \Lambda_1(A)$.
Given $A\in M_n$
and $t\in [0,2\pi)$, let $A(t) = e^{it}A+e^{-it}A^*$.  Consider the
compact convex sets
$$S(t) = \{\mu \in \Lambda_1(A): 
e^{it}\mu+ e^{-it}\overline{\mu}\le \lambda_k(A(t))\},
\qquad t \in [0, 2\pi).$$ 
By (\ref{lambda}),
$$\Lambda_k(A) = \bigcap_{t \in [0, 2\pi)} S(t).$$
By Helly's Theorem \cite[Theorem 24.9]{La}, 
it suffices to show that $S(t_1)\cap S(t_2)\cap S(t_3)\ne
\emptyset$ for all choices of $t_1, t_2, t_3$
with $0 \le t_1 < t_2 < t_3 < 2\pi$.

For $1\le j\le 3$, let $V_j$ be the subspace spanned by the
eigenvectors of $A(t_j)$ corresponding to the eigenvalues
$\lambda_k(A(t_j)), \dots , \lambda_n(A(t_j))$. Then $\dim V_j\ge
n-k+1$. Hence, we have
$$\begin{array}{rl}&\dim (V_1\cap V_2\cap V_3)\\
=&\dim (V_1\cap V_2)+\dim V_3-\dim((V_1\cap V_2)+ V_3)\\
=&\dim V_1+\dim V_2-\dim (V_1+ V_2)+\dim V_3-\dim((V_1\cap V_2)+
V_3)\\
\ge &3(n-k+1)-2n\\
=&n-3k+3\\
\ge &1.
\end{array}$$
Let $v$ be a unit (column) vector in $V_1\cap V_2\cap V_3$. Then
$\mu=v^*Av\in \Lambda_1(A)$ and for $1\le j\le 3$, we have
$$e^{it}\mu+ e^{-it}\overline{\mu}= v^*(A(t_j))v\le
\lambda_k(A(t_j)).$$ 
Hence, $\mu \in S(t_1)\cap S(t_2)\cap
S(t_3)$.\qed

The following answers a question in \cite{Cet0}.

\begin{corollary}
Let $A \in M_n$ with $n \ge 4$. Then $\Lambda_2(A) \ne \emptyset$.
\end{corollary}

Without additional information on $A \in M_n$, the bound on $n$ in
Theorem  1 is best possible as shown by the following
result.

\begin{theorem} \label{thm3}
Suppose $k$ is a positive integer such that  
$k \ge n/3+1$. There exists $A 
\in M_n$ such that
$\Lambda_k(A) = \emptyset$.
\end{theorem}

\it Proof. \rm
We first consider the case when $3k=n+3$.
Let $w = e^{i2\pi/3}$, and 
$$A = I_{k-1} \oplus  w I_{k-1} \oplus w^2 I_{k-1}.$$
Write $A = H+iG$ with $H = H^*$ and $G = G^*$.
Then $H = I_{k-1} \oplus (-1/2)I_{2k-2}$.
Thus, $\Lambda_k(H) = \{-1/2\}$; see also \cite[Theorem 2.4]{Cet1}. So,
$$\Lambda_k(A) \subseteq {\cal L} = \{z: \Re z = -1/2\}.$$
By rotation of $2\pi/3$ and $4\pi/3$, 
one can show that $\Lambda_k(A) \subseteq w{\cal L}$
and $\Lambda_k(A) \subseteq w^2 {\cal L}$.
So,
$$\Lambda_k(A) \subseteq {\cal L} \cap w{\cal L} \cap w^2 {\cal L} =
\emptyset.$$

\medskip
Now, suppose $3k>n+3$. Then we can consider a principal
submatrix $B \in M_n$ of the matrix $A \in M_{3k-3}$
constructed in the preceding paragraph.
Then $\Lambda_k(B) \subseteq \Lambda_k(A) = \emptyset$.
\qed

Note that we can perturb the example in the above proof to get a non-normal
matrix $A \in M_n$ such that $\Lambda_k(A) = \emptyset$
if $k\ge n/3+1$.  Also, Theorem \ref{thm3}
can be obtained from parts (1), (2), (3) of \cite[Theorem 4.7]{Cet0} 
and the fact that  $\Lambda_k(A)$ is 
a subset of 
$$\bigcap_{1 \le j_1 < \cdots < j_{n-k+1} \le n} 
\conv\{\lambda_{j_1},\dots, \lambda_{j_{n-k+1}}\}$$
if $A \in M_n$ is normal with eigenvalues $\lambda_1, \dots, \lambda_n$.

\medskip
Let $\BH$ be the algebra of bounded linear operator
acting on an infinite dimensional Hilbert space $\cH$.
One can extend the definition of $\Lambda_k(A)$
for a bounded linear operator $A \in \BH$ by
$$\Lambda_k(A) = \{\gamma \in \IC: X^*AX = \gamma I_k, \
X: \IC^k \rightarrow \cH,\ X^*X = I_k\}.$$
By Theorem 1, we have the following.

\begin{corollary}
Suppose $k$ is a positive integer
and $A \in \BH$ for an infinite dimensional Hilbert space
$\cH$.
Then
$$\Lambda_k(A) \ne \emptyset.$$
\end{corollary}

\bigskip\noindent
\bf Acknowledgment
\rm

\bigskip
We would like to thank the authors of \cite{Cet} and \cite{W} for
sending us their preprints. We also thank Professor John Holbrook
for some helpful comments, and the 
referee for a careful reading of the paper.

\end{document}